\documentclass[11pt,a4paper]{amsart}
\usepackage{amsmath,amssymb,amsthm}
\usepackage[margin=28mm]{geometry}
\usepackage{hyperref}
\theoremstyle{plain}
\newtheorem{theorem}{Theorem}[section]
\newtheorem{proposition}[theorem]{Proposition}
\newtheorem{lemma}[theorem]{Lemma}

\newtheorem{mainthm}{Theorem}

\theoremstyle{definition}

\newtheorem{remark}[theorem]{Remark}
\newtheorem{example}[theorem]{Example}

\newcommand{\R}{\mathbb R}

\newcommand{\D}{\mathbb D}
\newcommand{\HH}{\mathbb H}
\newcommand{\RHP}{\Pi}
\newcommand{\Lm}{\mathcal L_-}
\newcommand{\Lp}{\mathcal L_+}
\newcommand{\Mm}{\mathcal M_-}
\newcommand{\Mp}{\mathcal M_+}
\renewcommand{\Re}{\operatorname{Re}}
\renewcommand{\Im}{\operatorname{Im}}

\title[Kwong matrices and univalent maps]{Kwong matrices of operator monotone functions are infinitely divisible}

\newcommand{\AuthorName}{Tomoharu Suzuki}
\newcommand{\AuthorAddress}{Japan}
\author{\AuthorName}
\address{\AuthorAddress}
\email{quantale@gmail.com}
\date{November 2026}

\begin{document}

\begin{abstract}
For a function $f$ on $[0,\infty)$ and positive numbers $\lambda_1,\dots,\lambda_n$ the Kwong matrix is $[(f(\lambda_i)+f(\lambda_j))/(\lambda_i+\lambda_j)]$. Bhatia and Kosaki asked whether every function whose Loewner matrices are all infinitely divisible also has all its Kwong matrices infinitely divisible. We show more: the Kwong matrices of every nonnegative operator monotone function are infinitely divisible. Both ingredients are old: a nonconstant nonnegative operator monotone function extends to a univalent self-map of the right half-plane, which is essentially a theorem of Thale (1956) on Stieltjes transforms; and the kernel $(1-\Phi(z)\overline{\Phi(w)})/(1-z\bar w)$ is infinitely divisible for every univalent self-map $\Phi$ of the disc, a consequence of an inequality of Nehari (1953), stated in this form by Horn (1969). We include a short self-contained proof of both steps, and we show that the right half-plane is the largest domain on which all these functions are univalent. We also show that univalence is not necessary, with a one-parameter family whose threshold is sharp within the family, and that the Cayley conjugate of a real finite Blaschke product of degree at least two has Kwong matrices that are positive semidefinite but not all infinitely divisible. Finally, we restate the infinite divisibility of Kwong matrices as the contractivity of a composition operator on the Dirichlet space.
\end{abstract}

\maketitle

\section{Introduction}\label{sec:intro}

A real symmetric $n\times n$ matrix $A=[a_{ij}]$ with nonnegative entries is \emph{infinitely divisible} if the Hadamard power $A^{\circ r}=[a_{ij}^r]$ is positive semidefinite for every $r\ge0$ \cite{Ho69,BK} (we use the convention $0^0=1$, so that $A^{\circ0}$ is the all-ones matrix; Horn \cite[p.~272]{Ho69} takes $A^{\circ0}$ to be the incidence matrix instead, which makes no difference below). For a function $f$ on $[0,\infty)$ and positive numbers $\lambda_1,\dots,\lambda_n$ consider the matrices
\[
L_f=\Bigl[\frac{f(\lambda_i)-f(\lambda_j)}{\lambda_i-\lambda_j}\Bigr],\qquad K_f=\Bigl[\frac{f(\lambda_i)+f(\lambda_j)}{\lambda_i+\lambda_j}\Bigr],
\]
where the diagonal entries of $L_f$ are $f'(\lambda_i)$. We call $K_f$ the \emph{Kwong matrix} of $f$; such matrices are also called anti-Loewner matrices. Following Bhatia and Kosaki \cite[p.~49]{BK}, let $\mathcal F$ be the class of differentiable functions from $[0,\infty)$ into itself; throughout, we read this as continuous on $[0,\infty)$ and continuously differentiable on $(0,\infty)$, so that the powers $t^r$, $0<r<1$, belong to $\mathcal F$. Let $\mathcal L_\pm$ be the classes of $f\in\mathcal F$ for which all matrices $L_f$ (for $-$), respectively $K_f$ (for $+$), are positive semidefinite, and let $\mathcal M_\pm$ be the classes for which all these matrices are infinitely divisible. By Loewner's theorem \cite[Chap.~V]{BhMA} $\Lm$ is the class of nonnegative operator monotone functions on $[0,\infty)$ (here and below, operator monotone functions on $[0,\infty)$ are understood to be continuous at $0$, as members of $\mathcal F$ are; without this convention the indicator of $(0,\infty)$, which is operator monotone on $[0,\infty)$, would have to be admitted); by the theorem of Loewner, FitzGerald and Horn \cite{Lo66,Fi67,Ho67} (see \cite[p.~782]{Ho67} for the attribution) $\Mm$ consists of the constants and those $f\in\Lm$ whose analytic continuation maps the upper half-plane $\HH$ univalently into itself; and $\Lm\subset\Lp$ by Kwong \cite{Kw} (as quoted in \cite{BK}). Audenaert characterised $\Lp$ \cite{Au}. Ferrante \cite{Fe26} proves positive semidefiniteness of Cauchy-like matrices $[(r_i+r_j)/(\lambda_i+\lambda_j)]$ with $r_i$ depending only on the other $\lambda_j$. On $(0,1)$ instead of $(0,\infty)$, Kwong's theorem fails in general and holds for a subclass of operator monotone functions \cite{MST14}, and for the powers $t^r$, $r>1$, the inertia of the Kwong matrices is determined in \cite{BJ24}.

Bhatia and Kosaki asked \cite[p.~50]{BK}:

\smallskip
\noindent\emph{Question 1.} Is $\Mm\subset\Mp$?\quad \emph{Question 2.} Are there good characterisations of $\Lp$ and $\Mp$?
\smallskip

Known members of $\Mp$ are the powers $t^r$, $0\le r\le1$ (Bhatia--Kosaki \cite{BK}, after Kwong's proof \cite{Kw} that these $K_f$ are positive semidefinite; another proof, via Bapat's theorem, is in \cite{BS}) and the inverses of the functions $\beta t^3+\int t^3(\lambda+t^2)^{-1}d\nu(\lambda)$ (Hidaka--Sano \cite{HS}). We answer Question~1 in a stronger form.

\begin{mainthm}\label{thm:A}
Let $f$ be a nonnegative operator monotone function on $[0,\infty)$. Then for all positive numbers $\lambda_1,\dots,\lambda_n$ the matrix $K_f$ is infinitely divisible. In particular $\Mm\subset\Lm\subset\Mp$.
\end{mainthm}

Univalence is the mechanism, but it is not necessary:

\begin{mainthm}\label{thm:B}
For $0<a<1$ let $f_a(x)=\dfrac{(x+1)^2+a(x-1)^2}{(x+1)^2-a(x-1)^2}$. Then $f_a\in\Mp$ if and only if $a\le 1/\sqrt2$. Since $f_a(x)=f_a(1/x)$, the function $f_{1/\sqrt2}$ belongs to $\Mp$ but is not injective on $(0,\infty)$.
\end{mainthm}

\emph{On attribution.} Theorem~\ref{thm:A} is obtained by combining two results from the 1950s. The first is the univalence of nonconstant nonnegative operator monotone functions on the right half-plane (Lemma~\ref{lem:univ}), which is essentially a theorem of Thale on Stieltjes transforms \cite[Thm.~2.2]{Th56}. The second is the half-plane form (Theorem~\ref{thm:univ}) of the infinite divisibility of $(1-\Phi(z)\overline{\Phi(w)})/(1-z\bar w)$ for univalent self-maps $\Phi$ of the disc. This is an immediate consequence of Nehari's inequality \cite[(6)]{Ne}, which bounds the logarithm of this kernel, as a form in complex points with complex weights of sum zero, from below by the absolute value of a Grunsky-type form; Horn stated it in the form we use \cite[p.~285]{Ho69}, with a proof via the Grunsky inequalities in his Stanford thesis. We include short self-contained proofs of both steps (Lemma~\ref{lem:univ} and Section~\ref{sec:main}). What is new here is: (i) the observation that these two results together answer Question~1 of \cite{BK}, in the stronger form of Theorem~\ref{thm:A}; (ii) Section~\ref{sec:beyond}, which shows that univalence is not necessary (Theorem~\ref{thm:B}) and that Cayley conjugates of real finite Blaschke products of degree at least two have positive semidefinite Kwong matrices that are not all infinitely divisible; and (iii) Section~\ref{sec:remarks}, which restates membership in $\Mp$ as the contractivity of a composition operator on the Dirichlet space. In addition, Proposition~\ref{prop:max} shows that the right half-plane is the largest domain on which all these functions are univalent. The Cayley dictionary (Lemma~\ref{lem:dict}) is the Cayley transform of the half-plane Pick kernel; it is the tool for Sections~\ref{sec:beyond} and~\ref{sec:remarks}.

\section{Preliminaries}\label{sec:prelim}

\subsection{Infinite divisibility}
A real symmetric matrix $M$ is \emph{conditionally positive definite} if $\sum_{i,j}\xi_i\xi_jM_{ij}\ge0$ whenever $\xi\in\R^n$ and $\sum_i\xi_i=0$. A kernel $k$ is \emph{positive definite} if all matrices $[k(z_i,z_j)]$ are positive semidefinite; a positive semidefinite matrix with nonzero determinant is called \emph{nonsingular}.

\begin{lemma}[Loewner {\cite[p.~324]{Lo66}}, as quoted in {\cite[Cor.~1.6]{Ho69}}; see also {\cite[\S5.6]{PDM}}]\label{lem:cpd}
A real symmetric matrix $K$ with positive entries is infinitely divisible if and only if $[\log K_{ij}]$ is conditionally positive definite.
\end{lemma}

\begin{proof}
If $K$ is infinitely divisible and $\sum\xi_i=0$, then $\sum\xi_i\xi_j\frac{K_{ij}^r-1}{r}=\frac1r\sum\xi_i\xi_jK_{ij}^r\ge0$, and $r\downarrow0$ gives the claim. Conversely, let $\Lambda=[\log K_{ij}]$ and $P_{ij}=\Lambda_{ij}-\Lambda_{i1}-\Lambda_{1j}+\Lambda_{11}$. For $c\in\R^n$ the vector $\xi=c-(\sum_k c_k)e_1$ has $\sum\xi_i=0$ and $\xi^{\mathsf T}\Lambda\xi=c^{\mathsf T}Pc$, so $P\ge0$. With $d_i=\Lambda_{i1}-\frac12\Lambda_{11}$ we have $\Lambda_{ij}=P_{ij}+d_i+d_j$, hence $K^{\circ r}=D_r[\exp(rP_{ij})]D_r$ with $D_r=\operatorname{diag}(e^{rd_i})$, and $[\exp(rP_{ij})]=\sum_m\frac{r^m}{m!}P^{\circ m}\ge0$ by Schur's theorem.
\end{proof}

If $c_1,\dots,c_n>0$, then $K$ and $[c_ic_jK_{ij}]$ are infinitely divisible simultaneously (\emph{diagonal congruence}).

\subsection{The Cayley dictionary}
Let $\RHP=\{z:\Re z>0\}$, $\D$ the unit disc, $\psi(z)=\frac{z-1}{z+1}$ (a conformal map of $\RHP$ onto $\D$ with $\psi((0,\infty))=(-1,1)$), and for $f:(0,\infty)\to(0,\infty)$ put $\Phi=\psi\circ f\circ\psi^{-1}:(-1,1)\to(-1,1)$.

\begin{lemma}\label{lem:dict}
With $\zeta_x=\psi(x)$ and $\omega(x)=\frac{1+f(x)}{1+x}$,
\[
\frac{f(x)+f(y)}{x+y}=\omega(x)\omega(y)\,\frac{1-\Phi(\zeta_x)\Phi(\zeta_y)}{1-\zeta_x\zeta_y}\qquad(x,y>0).
\]
Consequently $K_f$ is infinitely divisible (resp.\ positive semidefinite) at $\lambda_1,\dots,\lambda_n$ if and only if the matrix $\bigl[\frac{1-\Phi(\zeta_i)\Phi(\zeta_j)}{1-\zeta_i\zeta_j}\bigr]$, $\zeta_i=\psi(\lambda_i)$, is.
\end{lemma}

\begin{proof}
For $a,b>-1$ one checks $1-\psi(a)\psi(b)=\frac{2(a+b)}{(1+a)(1+b)}$. Apply this to $(a,b)=(x,y)$ and to $(a,b)=(f(x),f(y))$, using $\psi(f(x))=\Phi(\zeta_x)$, and divide.
\end{proof}

The right-hand side of Lemma~\ref{lem:dict} is the Cayley transform of the half-plane Pick kernel $\frac{f(z)+\overline{f(w)}}{z+\bar w}$; we use it as a tool in Sections~\ref{sec:beyond} and~\ref{sec:remarks}. When $\Phi$ extends to a holomorphic self-map of $\D$, the kernel $k^\Phi(z,w)=\frac{1-\Phi(z)\overline{\Phi(w)}}{1-z\bar w}$ is the reproducing kernel of the de Branges--Rovnyak space $\mathcal H(\Phi)$; Lemma~\ref{lem:dict} says that $K_f$ is $k^\Phi$ at real points up to diagonal congruence. In particular all Kwong matrices of such an $f$ at positive points are positive semidefinite, and $f\in\Lp$ whenever $f\in\mathcal F$, that is, whenever $\Phi(\zeta)$ has a limit different from $1$ as $\zeta\downarrow-1$; this fails for $\Phi(\zeta)=\zeta^2$, which gives $f(x)=\frac{x^2+1}{2x}$. If $f$ is a nonconstant function in $\Lm$, then $\Phi=\psi\circ f\circ\psi^{-1}$, with $f$ continued to $\RHP$ by Lemma~\ref{lem:univ}, is a univalent self-map of $\D$. Conversely, let $f\in\Lp$ be not identically zero. Then $f>0$ on $(0,\infty)$: if $f(x_0)=0$ for some $x_0>0$, positive semidefiniteness of the $2\times2$ Kwong matrix at $x_0,y$ forces $f(y)=0$ for every $y>0$. Hence $\Phi$ is defined, and it extends to a holomorphic self-map of $\D$ with real Taylor coefficients: by Lemma~\ref{lem:dict} the matrices $[k^\Phi(\zeta_i,\zeta_j)]$ are positive semidefinite for all $\zeta_1,\dots,\zeta_n\in(-1,1)$, so by Pick's interpolation theorem \cite[Thm.~1.81]{AMY} $\Phi$ agrees on any finite subset of $(-1,1)$ with a holomorphic $\varphi:\D\to\overline{\D}$; by Montel's theorem a subsequence of such interpolants for an increasing sequence of finite sets with dense union converges locally uniformly to a holomorphic $\varphi:\D\to\overline{\D}$ which agrees with the continuous function $\Phi$ on $(-1,1)$, and $\varphi(\D)\subset\D$ by the maximum principle since $\Phi$ takes values in $(-1,1)$; its Taylor coefficients are real because it is real on $(-1,1)$.

\subsection{Univalence of operator monotone functions}
Every nonnegative operator monotone function $g$ on $[0,\infty)$ has the representation
\begin{equation}\label{eq:rep}
g(x)=\alpha+\beta x+\int_{(0,\infty)}\frac{x(1+t)}{x+t}\,d\mu(t),\qquad \alpha,\beta\ge0,
\end{equation}
with a finite positive measure $\mu$ \cite[Chapter V]{BhMA}, \cite{PDM}; $g$ is constant if and only if $\beta=0$ and $\mu=0$.

We call a function \emph{univalent} on a domain if it is holomorphic and injective there. The following lemma is essentially a theorem of Thale \cite[Thm.~2.2]{Th56}: if $\phi$ is nondecreasing and bounded on $[0,1]$ and not constant on $(0,1]$, then $F(x)=\int_{[0,1]}\frac{x}{1+x\tau}\,d\phi(\tau)$ is univalent on $\{\Re x>-1\}$. (Thale states this on the closed half-plane with $-1$ removed; there it can fail: if $\phi$ is concentrated at the two endpoints, with positive mass at each, then $F(x)=ax+bx/(1+x)$ has critical points on $\Re x=-1$, see the proof of Proposition~\ref{prop:max} with $t=1$; this is the only exception, because $\Im F'(x)=-2\Im x\int\frac{\tau(1+\tau\Re x)}{\lvert1+x\tau\rvert^4}\,d\phi(\tau)$, whose integrand is nonnegative for $\Re x\ge-1$ and vanishes on $\Re x=-1$ only at $\tau\in\{0,1\}$; so if $\phi$ charges $(0,1)$, then $\Im F'<0$ on the convex set $\{\Re x\ge-1,\ \Im x>0\}$, and integrating $F'$ along segments shows that $F$ is injective there (and likewise in the lower half), while $\Im F$ has the sign of $\Im x$ and $F'>0$ on $(-1,\infty)$. His argument gives the open half-plane, which is all we use.) Indeed, with $\tau=1/t$ we have $\frac{x(1+t)}{x+t}=(1+\tau)\frac{x}{1+x\tau}$, and the term $\beta x$ corresponds to a point mass at $\tau=0$. If $\mu$ is supported in $[\varepsilon,\infty)$, the substitution $x=\varepsilon y$, $\tau=\sigma/\varepsilon$ brings $g-\alpha$ to Thale's form in $y$ (the excluded case being $g-\alpha=\beta x$), so $g$ is univalent on $\{\Re x>-\varepsilon\}\supset\RHP$. In general, if $\mu\ne0$ (otherwise $g$ is affine), $g$ is the locally uniform limit on $\RHP$ of the functions obtained by restricting $\mu$ to $[\varepsilon,\infty)$, and Hurwitz's theorem applies because $g$ is nonconstant. We give a direct proof, in the spirit of Thale's.

\begin{lemma}\label{lem:univ}
Let $g$ be a nonconstant nonnegative operator monotone function on $[0,\infty)$. Then \eqref{eq:rep} defines a holomorphic extension of $g$ to $\RHP$ with $g(\RHP)\subset\RHP$, and $g$ is univalent on $\RHP$.
\end{lemma}

\begin{proof}
For $z\in\RHP$ and $t>0$ we have $\lvert z+t\rvert\ge\Re z+t$, so the integrand in \eqref{eq:rep} is bounded by $\lvert z\rvert\max(1,1/\Re z)$, locally uniformly on $\RHP$; dominated convergence and Morera's theorem give holomorphy. Since $\Re\frac{z}{z+t}=\frac{\lvert z\rvert^2+t\Re z}{\lvert z+t\rvert^2}>0$ and $\beta>0$ or $\mu\ne0$, we get $\Re g(z)>0$.

For injectivity write $\frac{x(1+t)}{x+t}=(1+t)-\frac{t(1+t)}{x+t}$. For $z\ne w$ in $\RHP$,
\[
\frac{g(z)-g(w)}{z-w}=\beta+\int_{(0,\infty)}\frac{t(1+t)}{(z+t)(w+t)}\,d\mu(t),
\]
the integrand being bounded in $t$. Let $\theta_1=\arg(z+t)$ and $\theta_2=\arg(w+t)$; both lie in $(-\frac\pi2,\frac\pi2)$, and the integrand has argument $-(\theta_1+\theta_2)$.
If $\Im z,\Im w\ge0$, then $\theta_1+\theta_2\in[0,\pi)$, so the integrand lies in the closed lower half-plane, and it is a positive real number wherever its imaginary part vanishes. Hence the right-hand side has nonpositive imaginary part; if this imaginary part is $0$, then the integrand is a positive real number for $\mu$-almost every $t$, and the right-hand side is positive because $\beta>0$ or $\mu\ne0$. If $\Im z,\Im w\le0$, use $g(\bar z)=\overline{g(z)}$. If $\Im z>0>\Im w$, then $\theta_1+\theta_2\in(-\frac\pi2,\frac\pi2)$ and the integrand has positive real part. These cases cover all pairs $z\ne w$ in $\RHP$ (a pair containing a real point falls under one of the first two), and in every case the right-hand side is nonzero.
\end{proof}

\begin{proposition}\label{prop:max}
Let $D$ be a domain with $D\cap(0,\infty)\ne\emptyset$. Every nonconstant nonnegative operator monotone function extends to a univalent holomorphic function on $D$ if and only if $D\subset\RHP$. In particular, for no $\theta>\frac\pi2$ is every such function univalent on the sector $\{\lvert\arg z\rvert<\theta\}$.
\end{proposition}

\begin{proof}
If $D\subset\RHP$, restrict Lemma~\ref{lem:univ}. Conversely, if $D\not\subset\RHP$, then $D$, being open, contains a point $z_0=-s+iy$ with $s>0$ and $y\ne0$. Let $t=s$, $c=y^2/s$ and $g(x)=x+\frac{cx}{x+t}$, which is of the form \eqref{eq:rep} with $\alpha=0$, $\beta=1$, $\mu=\frac{c}{1+t}\delta_t$. Since $g'(z)=1+\frac{ct}{(z+t)^2}$, the critical points of $g$ are $-t\pm i\sqrt{ct}$, and $g'(z_0)=0$. A holomorphic extension $G$ of $g$ to $D$ coincides with the rational function $g$ on the connected set $D\setminus\{-t\}$ by the identity theorem. If $-t\in D$ this contradicts the pole of $g$; otherwise $G'(z_0)=0$, so $G$ is not injective near $z_0$. For the sector take $t=1$ and $c>\tan^2\theta$.
\end{proof}

Thale \cite[Thm.~2.2]{Th56} shows that no region containing $\{\Re x\ge-1\}\setminus\{-1\}$ and a further point is a region of univalence for the whole class considered there. Proposition~\ref{prop:max} is the corresponding statement for the class \eqref{eq:rep}, for every domain meeting $(0,\infty)$ and not contained in $\RHP$, with an explicit function for each excluded point; the example $x+cx/(x+1)$ also shows that the line $\Re x=-1$ itself must be excluded from Thale's region when $\phi$ is concentrated at the two endpoints.

\section{Univalent self-maps and the proof of Theorem~\ref{thm:A}}\label{sec:main}

Let $A^2(\RHP)$ be the Bergman space of holomorphic $u$ on $\RHP$ with $\lVert u\rVert^2=\int_\RHP\lvert u\rvert^2dA<\infty$. Its reproducing kernel is $\frac{1}{\pi(z+\bar w)^2}$; in particular $u(s)=\frac1\pi\int_\RHP u(z)\,\overline{(z+s)^{-2}}\,dA(z)$ for $s>0$. For $a=(a_1,\dots,a_n)\in(0,\infty)^n$ and real $\xi$ with $\sum\xi_i=0$ put
\[
F_a(z)=\sum_i\frac{\xi_i}{z+a_i},\qquad Q(a)=-\sum_{i,j}\xi_i\xi_j\log(a_i+a_j).
\]
Since $\sum\xi_i=0$, $F_a(z)=O(\lvert z\rvert^{-2})$ at infinity, so $F_a\in A^2(\RHP)$.

\begin{lemma}\label{lem:repr}
If $u\in A^2(\RHP)$ and $U$ is holomorphic on $\RHP$ with $U'=u$, then $\sum_i\xi_iU(a_i)=-\frac1\pi\langle u,F_a\rangle$. In particular $\lVert F_a\rVert^2=\pi Q(a)\ge0$.
\end{lemma}

\begin{proof}
Fix $b>0$. Since $\sum\xi_i=0$, $\sum_i\xi_iU(a_i)=\sum_i\xi_i\int_b^{a_i}u(s)\,ds$. Insert the reproducing formula and use Fubini (the integrand is integrable on $[\min,\max]\times\RHP$): $\int_b^{a}\frac{ds}{(\bar z+s)^2}=\frac1{\bar z+b}-\frac1{\bar z+a}$, and the terms with $b$ cancel because $\sum\xi_i=0$. For the second claim take $u=F_a$ and $U(z)=\sum_j\xi_j\log(z+a_j)$ (principal branch).
\end{proof}

\begin{theorem}[Nehari {\cite[(6)]{Ne}}, Horn {\cite[p.~285]{Ho69}}; half-plane form]\label{thm:univ}
Let $f$ be a univalent map of $\RHP$ into itself which is positive on $(0,\infty)$. Then for all positive $x_1,\dots,x_n$ the matrix $\bigl[\log\frac{f(x_i)+f(x_j)}{x_i+x_j}\bigr]$ is conditionally positive definite, and $K_f$ is infinitely divisible.
\end{theorem}

\begin{proof}
Let $a_i=f(x_i)$ and let $\xi$ be real with $\sum\xi_i=0$; we must show $Q(a)\le Q(x)$. Put $u=(F_a\circ f)\,f'$. Since $f$ is injective, the change of variables formula gives
$\lVert u\rVert^2=\int_\RHP\lvert F_a(f(z))\rvert^2\lvert f'(z)\rvert^2dA(z)=\int_{f(\RHP)}\lvert F_a\rvert^2dA\le\lVert F_a\rVert^2$, so $u\in A^2(\RHP)$, and $u=U'$ with $U(z)=\sum_j\xi_j\log(f(z)+a_j)$. By Lemma~\ref{lem:repr} applied at the points $x_i$,
\[
\pi Q(a)=-\pi\sum_i\xi_iU(x_i)=\langle u,F_x\rangle\le\lVert u\rVert\,\lVert F_x\rVert\le\lVert F_a\rVert\,\lVert F_x\rVert=\pi\sqrt{Q(a)Q(x)}.
\]
As $Q(a)\ge0$, this gives $Q(a)\le Q(x)$, i.e.\ $\sum\xi_i\xi_j\log\frac{a_i+a_j}{x_i+x_j}=Q(x)-Q(a)\ge0$. Lemma~\ref{lem:cpd} finishes the proof.
\end{proof}

\begin{remark}
Injectivity is used only in the area estimate; it is sufficient but not necessary (Theorem~\ref{thm:B}). For $n=2$ no injectivity is needed: the inequality is the Schwarz--Pick lemma for $\RHP$, i.e.\ $\rho(f(x),f(y))\le\rho(x,y)$ with $\rho(z,w)=\lvert\frac{z-w}{z+\bar w}\rvert$. The area argument is in the spirit of Nehari's derivation of \cite[(6)]{Ne} from the minimum property of the Dirichlet integral.
\end{remark}

\begin{proof}[Proof of Theorem~\ref{thm:A}]
If $f\equiv0$ all matrices vanish and are infinitely divisible, with the convention $0^0=1$. If $f\equiv c>0$, then $[\log K_{ij}]=[\log 2c-\log(\lambda_i+\lambda_j)]$ is conditionally positive definite by Lemma~\ref{lem:repr} ($Q(\lambda)\ge0$), so $K_f$ is infinitely divisible by Lemma~\ref{lem:cpd}. If $f$ is nonconstant, Lemma~\ref{lem:univ} shows that $f$ satisfies the hypotheses of Theorem~\ref{thm:univ}. Finally $\Mm\subset\Lm$ trivially.
\end{proof}

\begin{remark}[The Hidaka--Sano functions]\label{rem:HS}
Let $\phi(t)=\beta t^3+\int_{(0,\infty)}\frac{t^3}{\lambda+t^2}d\nu(\lambda)$, nonzero, and $g=\phi^{-1}$ on $(0,\infty)$. Hidaka and Sano \cite{HS} proved that $K_\phi$ is conditionally negative definite and deduced that $K_g$ is infinitely divisible. (Since $K_g(\phi(s),\phi(t))=1/K_\phi(s,t)$, the deduction can be made by Bapat's theorem that the entrywise reciprocal of a conditionally negative definite matrix with positive entries is infinitely divisible, see \cite{BS}; the kernel form of this mechanism, and its use for inverse functions, go back to Horn \cite[Cor.~4, Thm.~5]{Ho71}.) The infinite divisibility of $K_g$ also follows from Theorem~\ref{thm:univ}, although $g$ need not be operator monotone. Indeed, $\phi(t)=t\,h(t^2)$ with $h(s)=\beta s+\int\frac{s}{\lambda+s}d\nu$ nonnegative operator monotone. Let $p(s)=g(\sqrt s)^2$; it is the inverse of $q(y)=y\,h(y)^2$. By Uchiyama's product formula $\mathcal P\cdot\mathcal P^{-1}\subset\mathcal P^{-1}$ \cite[Lemma~2.4]{Uc06} (with $\mathcal P$ the nonnegative operator monotone functions on $[0,\infty)$ and $\mathcal P^{-1}$ the continuous strictly increasing functions $k$ on $[0,\infty)$ with $k([0,\infty))=[0,\infty)$ and $k^{-1}$ operator monotone on $[0,\infty)$), applied twice to $y\in\mathcal P^{-1}$ and $h\in\mathcal P$, we get $q\in\mathcal P^{-1}$, that is, $p=q^{-1}$ is operator monotone on $[0,\infty)$. The continuation of $p$ to $\HH$ maps $\HH$ into $\HH$, and $q\circ p=\mathrm{id}$ on $\HH$ by the identity theorem, so $p$ is univalent on $\HH$. Since $0\le\arg p(w)\le\arg w$ on $\HH$ (each term of \eqref{eq:rep} has argument in $[0,\arg w]$), $z\mapsto\sqrt{p(z^2)}$ maps the first quadrant injectively into itself and agrees with $g$ on $(0,\infty)$; by the reflection principle $g$ extends to a univalent self-map of $\RHP$.
\end{remark}

\section{Beyond univalence}\label{sec:beyond}

\begin{proof}[Proof of Theorem~\ref{thm:B}]
The denominator of $f_a$ is positive on $[0,\infty)$, so $f_a$ maps $[0,\infty)$ smoothly into $(0,\infty)$. Since $f_a=\psi^{-1}\circ(a\psi^2)$, its Cayley conjugate is $\Phi(\zeta)=a\zeta^2$. By Lemma~\ref{lem:dict}, with $\tau=\zeta_x\zeta_y\in(-1,1)$,
\[
\log K_{f_a}(x,y)=\log\omega(x)+\log\omega(y)+\log\frac{1-a^2\tau^2}{1-\tau}=\log\omega(x)+\log\omega(y)+\sum_{m\ge1}b_m\tau^m,
\]
with $b_{2k-1}=\frac1{2k-1}$ and $b_{2k}=\frac{1-2a^{2k}}{2k}$.

If $a^2\le\frac12$, all $b_m\ge0$. For $r>0$, $\exp(r\sum_m b_m\tau^m)$ is a power series in $\tau$ with nonnegative coefficients, convergent for $\lvert\tau\rvert<1$; since each $[\zeta_i^m\zeta_j^m]$ is positive semidefinite, so is $[\exp(r\sum_mb_m(\zeta_i\zeta_j)^m)]$, and $K_{f_a}^{\circ r}$ is diagonally congruent to it.

If $a^2>\frac12$, take $\eta=(-1,0,1)$, $\xi=(1,-2,1)$ and, for $0<\varepsilon<1$, the points $x_i=\psi^{-1}(\varepsilon\eta_i)>0$. Then $\sum\xi_i=\sum\xi_i\eta_i=0$, $\sum\xi_i\eta_i^2=2$, and
\[
\sum_{i,j}\xi_i\xi_j\log K_{f_a}(x_i,x_j)=\sum_mb_m\varepsilon^{2m}\Bigl(\sum_i\xi_i\eta_i^m\Bigr)^2=4\sum_{k\ge1}b_{2k}\varepsilon^{4k}=2\log\frac{(1-a^2\varepsilon^4)^2}{1-\varepsilon^4},
\]
because $\sum_i\xi_i\eta_i^m$ equals $0$ for odd $m$ and $2$ for even $m$. This is negative if and only if $a^4\varepsilon^4<2a^2-1$, which holds for all small $\varepsilon>0$ since $a^2>\frac12$. By Lemma~\ref{lem:cpd}, $K_{f_a}$ is not infinitely divisible. Finally $\psi(1/x)=-\psi(x)$ gives $f_a(1/x)=f_a(x)$.
\end{proof}

\begin{proposition}\label{prop:blaschke}
Let $\Phi$ be a finite Blaschke product of degree $d\ge2$ with real coefficients and $f=\psi^{-1}\circ\Phi\circ\psi$ on $(0,\infty)$. Then all Kwong matrices of $f$ at positive points are positive semidefinite, and $K_f$ is not infinitely divisible at any $d+1$ distinct positive points. If moreover $\Phi(-1)=-1$, then $f$ extends continuously to $[0,\infty)$ with $f(0)=0$, and $f\in\Lp\setminus\Mp$.
\end{proposition}

\begin{proof}
$f$ maps $(0,\infty)$ into itself, and positivity follows from Lemma~\ref{lem:dict} since $k^\Phi$ is a positive definite kernel. The space $\mathcal H(\Phi)$ is the model space $H^2\ominus\Phi H^2$, of dimension $d$, so $k^\Phi$ and hence $K_f$ have rank at most $d$. At two distinct points the $2\times2$ minor of $k^\Phi$ is positive: it is positive if and only if $\rho(\Phi(z),\Phi(w))<\rho(z,w)$ for the pseudo-hyperbolic distance $\rho(z,w)=\bigl\lvert\frac{z-w}{1-z\bar w}\bigr\rvert$ of $\D$, which holds by the Schwarz--Pick lemma because $\Phi$ is not an automorphism. By Farjot's theorem \cite[Thm.~2.1]{Fa}, an infinitely divisible positive semidefinite matrix of order $n>2$ all of whose $2\times2$ principal minors are positive is nonsingular. At $d+1\ge3$ distinct points this would force rank $d+1$, a contradiction. Finally, if $\Phi(-1)=-1$, then $f(x)=\psi^{-1}(\Phi(\psi(x)))\to\psi^{-1}(-1)=0$ as $x\downarrow0$, so $f\in\mathcal F$, and the two assertions just proved give $f\in\Lp\setminus\Mp$.
\end{proof}

\begin{example}\label{ex:p}
For $p(x)=\frac{x}{1+x^2}$ one finds $\Phi(\zeta)=-\frac{1+3\zeta^2}{3+\zeta^2}$, a Blaschke product of degree $2$ with zeros $\pm i/\sqrt3$ and $\Phi(-1)=-1$; equivalently $K_p(x,y)=\frac{1+xy}{(1+x^2)(1+y^2)}$ is diagonally congruent to $[1+x_ix_j]$. Thus $p\in\Lp\setminus\Mp$ (this also follows from Jain \cite[Thm.~1.1]{Ja17}). An explicit certificate: at $(\frac12,1,2)$ with $\xi=(1,-2,1)$,
\[
\sum_{i,j}\xi_i\xi_j\log K_p(x_i,x_j)=\log\frac{(4/5)(1/2)^4(1/5)(8/25)^2}{(3/5)^4(3/10)^4}=\log\frac{6400}{6561}<0.
\]
\end{example}

\begin{remark}[Strict inclusions]\label{rem:chain}
Let $\mathcal U$ denote the functions in $\mathcal F$ which extend to univalent self-maps of $\RHP$. Then $\Lm\setminus\{\text{constants}\}\subset\mathcal U\subset\Mp$ (Lemma~\ref{lem:univ}, Theorem~\ref{thm:univ}), and all inclusions in $\Mm\subset\Lm\subset\Mp\subset\Lp$ are strict, as are $\Lm\setminus\{\text{constants}\}\subset\mathcal U\subset\Mp$. Witnesses:
$\frac{x}{1+x}+\frac{x}{3(3+x)}\in\Lm\setminus\Mm$ (its derivative $(1+z)^{-2}+(3+z)^{-2}$ vanishes at $-2+i\in\HH$, so it is not univalent on $\HH$; numerically, at $(\frac1{20},2,\frac{12}5,20)$ with $\xi=(1,-1,-1,1)$ the form $\sum\xi_i\xi_j\log(L_f)_{ij}$ equals $-0.00267\ldots$);
$\frac1{1+x}\in\mathcal U\setminus\Lm$;
$f_{1/\sqrt2}\in\Mp\setminus\mathcal U$ (Theorem~\ref{thm:B});
$\frac{x}{1+x^2}\in\Lp\setminus\Mp$ (Example~\ref{ex:p}). Note that a function in $\mathcal U$ is automatically positive on $(0,\infty)$, since $\RHP\cap\R=(0,\infty)$.
\end{remark}

\section{Remarks on Question 2}\label{sec:remarks}
Since $K_{cf}=cK_f$ for $c>0$, for $f\ne0$ membership in $\Mp$ is unchanged by normalising $f(1)=1$, i.e.\ $\Phi(0)=0$. Let $\mathcal D_0$ be the Dirichlet space of holomorphic $h$ on $\D$ with $h(0)=0$ and $\lVert h\rVert^2=\frac1\pi\int_\D\lvert h'\rvert^2dA<\infty$; its reproducing kernel is $s(z,w)=-\log(1-z\bar w)=\sum_{m\ge1}\frac{z^m\bar w^m}{m}$.

\begin{proposition}\label{prop:dirichlet}
Let $\Phi$ be a holomorphic self-map of $\D$ with real Taylor coefficients and $\Phi(0)=0$, and let $\ell(z,w)=s(z,w)-s(\Phi(z),\Phi(w))$. Since $1-z\bar w$ and $1-\Phi(z)\overline{\Phi(w)}$ lie in $\RHP$, $\ell$ is the principal branch of $\log\frac{1-\Phi(z)\overline{\Phi(w)}}{1-z\bar w}$, holomorphic in $(z,\bar w)$ on $\D\times\D$. The following are equivalent:
\begin{enumerate}
\item[(a)] $[\ell(\zeta_i,\zeta_j)]$ is conditionally positive definite for all $\zeta_1,\dots,\zeta_n\in(-1,1)$;
\item[(b)] $\ell$ is a positive definite kernel on $\D$;
\item[(c)] the composition operator $C_\Phi h=h\circ\Phi$ is a contraction on $\mathcal D_0$.
\end{enumerate}
\end{proposition}

\begin{proof}
(b)$\Rightarrow$(a) is clear. (a)$\Rightarrow$(b): since $\ell(\cdot,0)=0$, adjoining the point $0$ with weight $-\sum_i c_i$ shows that $\sum_{i,j}c_ic_j\ell(\zeta_i,\zeta_j)\ge0$ for all real $c$, so $\ell$ is a positive definite kernel on $(-1,1)$. Now $\ell(z,w)=\sum_{j,k\ge1}a_{jk}z^j\bar w^k$ with real $a_{jk}$, converging absolutely on $\D\times\D$ (by Cauchy's estimates for the coefficients of $\Phi^m$, since $\max_{\lvert z\rvert\le r}\lvert\Phi\rvert<1$ for $r<1$). For $k\ge1$ and small $\delta>0$ let $c^{(k)}_p=\frac{(-1)^{k-p}}{k!\,\delta^k}\binom kp$ for $0\le p\le k$, so that $\sum_pc^{(k)}_pg(p\delta)\to g^{(k)}(0)/k!$ for every smooth $g$. Given real $\gamma_1,\dots,\gamma_N$, apply positive definiteness on $(-1,1)$ to the points $p\delta$, $0\le p\le N$, with weights $c_p=\sum_k\gamma_kc^{(k)}_p$, and let $\delta\to0$; since $\ell$ is real analytic on $(-1,1)^2$, the double divided differences converge and we get $\sum_{j,k=1}^N\gamma_j\gamma_k\,\frac{\partial_x^j\partial_y^k\ell(0,0)}{j!\,k!}=\sum_{j,k=1}^N\gamma_j\gamma_ka_{jk}\ge0$. Thus every finite section of $[a_{jk}]$ is positive semidefinite. Each partial sum $\sum_{j,k\le N}a_{jk}z^j\bar w^k$ is therefore a Gram kernel on $\D$, and so $\ell$, their pointwise limit, is a positive definite kernel on $\D$.
(b)$\Leftrightarrow$(c): on the span of the kernels $s_w=s(\cdot,w)$ define $Ts_w=s_{\Phi(w)}$. Then
$\lVert\sum c_is_{w_i}\rVert^2-\lVert\sum c_is_{\Phi(w_i)}\rVert^2=\sum_{i,j}\bar c_ic_j\,\ell(w_i,w_j)$,
so (b) says exactly that $T$ is well defined and contractive; it then extends to a contraction on $\mathcal D_0$, and $T^{\ast}h(w)=\langle h,s_{\Phi(w)}\rangle=h(\Phi(w))$, i.e.\ $T^{\ast}=C_\Phi$. Conversely if $C_\Phi$ is a contraction, then $C_\Phi^{\ast}s_w=s_{\Phi(w)}$ and the displayed identity gives (b).
\end{proof}

The equivalence (b)$\Leftrightarrow$(c) is the case $c=1$ of the general fact that $\lVert C_\Phi\rVert$ is the least $c\ge0$ for which $c^2s(z,w)-s(\Phi(z),\Phi(w))$ is a positive definite kernel (see \cite[Theorems~2.9 and~2.17]{SS16}, and \cite{Ju07} for this approach on Hardy and Bergman spaces); the new part is (a)$\Rightarrow$(b).
By Lemmas~\ref{lem:dict} and~\ref{lem:cpd}, condition (a) holds if and only if all Kwong matrices of $f=\psi^{-1}\circ\Phi\circ\psi$ at positive points are infinitely divisible; for $f\in\mathcal F$, that is, when $\Phi(\zeta)$ has a limit different from $1$ as $\zeta\downarrow-1$, this means $f\in\Mp$.

Since the Cayley conjugate of every nonzero $f\in\Lp$ extends to a self-map of $\D$ with real Taylor coefficients (Section~\ref{sec:prelim}), Question~2 for $\Mp$ asks which $\Phi$ give contractive composition operators on the Dirichlet space. Univalence suffices (by the area formula, as in the proof of Theorem~\ref{thm:univ}; this is standard for composition operators on the Dirichlet space); it is not necessary (Theorem~\ref{thm:B}); and finite Blaschke products of degree $\ge2$ are excluded (Proposition~\ref{prop:blaschke}). To our knowledge, no function-theoretic characterisation of the contractive composition operators on $\mathcal D_0$ is known; for the isometries, see Mart\'{\i}n and Vukoti\'c \cite{MV06}. The norm of $C_\Phi$ gives quick checks: $C_{a\zeta^2}$ maps $z^m$ to $a^mz^{2m}$ and $\lVert z^m\rVert^2=m$, so $\lVert C_{a\zeta^2}\rVert^2=\sup_m2a^{2m}=2a^2$, which recovers the threshold of Theorem~\ref{thm:B} (note $f_a(1)=1$); and if $\Phi$ is a finite Blaschke product of degree $d\ge2$ with $\Phi(0)=0$, then $\lVert C_\Phi z\rVert^2=\lVert\Phi\rVert^2=d>1$ by the area formula, so (c) fails, in accordance with Proposition~\ref{prop:blaschke}.

\section*{Declaration of competing interest}
The author declares that they have no known competing financial interests or personal relationships that could have appeared to influence the work reported in this paper.

\section*{Data availability}
No data was used for the research described in the article.

\section*{Funding}
This research did not receive any specific grant from funding agencies in the public, commercial, or not-for-profit sectors.

\section*{Declaration of generative AI and AI-assisted technologies in the manuscript preparation process}
During the preparation of this work the author used Claude (Anthropic) in order to search for and organise related literature, to write code for numerical checks (on which none of the proofs depends), and to improve the language of the manuscript. All references were checked by the author against the original sources. After using this tool, the author reviewed and edited the content as needed and takes full responsibility for the content of the published article.

\end{document}